\documentclass{amsart}
\begin{document}
\author {N. Mohan Kumar}
\address{Department of Mathematics\\
Washington University in St. Louis}
\email{ kumar@math.wustl.edu}
\thanks{Work partially supported by NSF grants}
\title{A note on cancellation of reflexive modules}
\newtheorem{thm}{Theorem}
\newtheorem{lem}{Lemma}
\newtheorem{cor}{Corollary}
\maketitle
\section{Introduction}
By the Quillen-Suslin theorem \cite{Q,sus4}, we know that
projective modules over a polynomial ring over a field are free. One
way of saying this is, that if two projective modules of the same rank
are stably
isomorphic, then they are isomorphic. That, projective modules of
given rank over polynomial rings are stably isomorphic was well known
at the time Quillen and Suslin proved their theorems
and this result is usually attributed to Grothendieck. This result can
also be deduced from Hilbert Syzygy Theorem. This note tries to answer
whether a similar 
cancellation occurs for reflexive modules. The question was
specifically raised by M. P. Murthy. I thank him for raising this
question and for the innumerable discussions which ensued. In this
note we show that, in general, reflexive modules are not
cancellative, without further assumptions. These assumptions under
which cancellation does take place are explained in the theorem in the
next section, Theorem \ref{atheorem}. 

\section{A case where cancellation is true}
As usual, we say that a module $M$ over a ring $A$ is
cancellative, if for some finitely generated free module $F$ over $A$
and a module $N$, if $F\oplus M\cong F\oplus N$, then $M\cong N$.
 The main
theorem we will  prove is the following:
\begin{thm}\label{atheorem}
Let $R$ be an affine domain over an algebraically closed field
of characteristic zero,
$M$ a reflexive module over $R$ of finite homological dimension such
that $M$ is locally free outside 
a finite set of closed points of $R$. Further assume that
${\rm rank}\, M\geq \dim R$. Then $M$ is
cancellative. 
\end{thm}

We start with some lemmas. We will not state the most general
versions,  but just what we need. 

\begin{lem}\label{lem1}
Let $A$ be an affine algebra over an infinite field $k$ of dimension
$n$ and $Q$ a projective module (of constant rank) over $A$. Let
$a\in A$ and let $I\subset A$ be an ideal of height $n$ . Then there
exists an element $f=\lambda+a+x$, with 
$x\in I$ and $0\neq \lambda\in k$ such that $Q_f$ is free over $A_f$.
\end{lem}

\medskip
\noindent{\bf Proof:} Since $I$ is of height $n$, $Q$ is free when we
semi-localise at the finitely many maximal ideals containing $I$. Thus
we may find an $s\in
A$, comaximal with $I$ such that $Q_s$ is $A_s$-free. Since the base
field is infinite, for a general $0\neq \lambda\in k$,
$a+\lambda$ is comaximal with $I$. Since  $s$ is a unit
modulo $I$, we can find a $t\in A$ such that $st\equiv a+\lambda\pmod
I$. Notice that $Q_{st}$ is clearly free over $A_{st}$.
So, letting $f=st$, we have,
$f=\lambda+a+x$ with $x\in I$.

\begin{cor}\label{coro1}
Notation being the same as in the lemma, there exists an $x\in I$ such
that either $a+x$ is a unit in $A$ or $Q/(a+x)Q$ is a free module over
$A/(a+x)A$. 
\end{cor}

\medskip
\noindent{\bf Proof:} From the lemma, we have $f-\lambda=a+x$ with
$0\neq \lambda\in k$ and $Q_f$ is $A_f$ free. But, then
$Q/(f-\lambda)Q$ is clearly free over $A/(f-\lambda)A$, unless
$f-\lambda$ is a unit. 

In the following two lemmas, $A$
will be an integral domain, $M$ a finitely generated torsion-free
$A$-module, $0\neq a\in A, m\in M$ such that $M/aM$ is a projective
$A$-module of rank equal to the rank of $M=n$ over $A/aA$  and the
image of $m$ in $M/aM$ is a unimodular element of $M/aM$. The
conditions imply (by an easy local checking), that the map
$A\stackrel{(a,m)}{\longrightarrow} A\oplus M$ is a split inclusion.

\begin{lem}\label{lem2}
Let 
$m'\in M$ and $d\in\mathbb{N}$. Then, 
$$A\oplus M/(a^d,m)A\cong A\oplus M/(a^d, m+am')A.$$
\end{lem}

The proof is essentially the same as  in \cite[lemma 2]{MK97} and is
just a slight modification of a result of Suslin \cite{sus2}. 

\medskip
\noindent{\bf Proof:} Consider $B=A[t]$, polynomial ring in one variable
over $A$ and let us consider the module, $N=B\oplus M[t]/(a^d,
m+atm')$. We will show that $N$ is an extended module and then putting
$t=0,1$, we would be done.

The hypotheses on $a,m$ imply  that $B\oplus N\cong B\oplus
M$. So, if $n=1$, we are done by taking determinants. So, we
will further assume that  $n\geq 2$.

To check that $N$ is extended, by Quillen's theorem \cite{Q}, we
need to show this only locally. At maximal ideals not containing $a$,
clearly $N\cong M[t]$. At maximal ideals containg $a$, by
hypothesis $M$ is free with $m$ as part of a basis. So choose a basis,
$m=m_1, m_2,\ldots, m_n$. Write $m'=\sum c_im_i$. Then with respect to
this basis, the vector $(a^d, m+atm')$ is,
$(a^d, 1+ac_1t, ac_2t,\ldots, ac_nt).$ Since $a^d, 1+ac_1t$ generate
the unit ideal, we can change $ac_2t$ (which exists, since $n\geq 2$)
by elementary transformation to 1 and thus $N$ is free at such a
maximal ideal.

Next we prove a crucial lemma, which is essentially a slight
generalisation of a theorem of Suslin \cite{sus1}.
\begin{lem}\label{lem3}
Notation being as above,  assume further that $M/aM$ is a free module
and the image of $m$ in $M/aM$ is part of a free basis of
$M/aM$. Then, $$A\oplus M/(a^n,m)\cong M.$$
\end{lem}

\medskip
\noindent{\bf Proof:} Since $M/aM$ is free over $A/aA$ with image of
$m$ as a part of a free basis, we may choose $m=m_1,m_2,\ldots, m_n\in
M$ such that their images in $M/aM$ form a free basis. Consider $B=A[t]$
and maps, $\phi(t):B^n\to B^n$ and $\psi:B^n\to M$ given as follows. The
map $\psi$ is just sending a basis $\{e_i\}$ of $B^n$ to the
$m_i$'s. The map 
$\phi(t)$ is given by the $n\times n$ matrix, which has $a$ for its
diagonal entries and $t$ on the subdiagonal, with zero
elsewhere. That is, $\phi(t)(e_1)=ae_1$ and $\phi(t)(e_i)= te_{i-1}+ae_i$
for $i>1$. Consider $N=B^n\oplus M[t]/K$ where $K$ is the image of
$B^n$ under the map $(\phi(t),\psi)$. 

First, I claim that $N$ is extended. Again, by Quillen's theorem {\em
loc. cit}, suffices to do this locally on $A$. At a maximal ideal not
containg $a$, $\phi(t)$ is an isomorphism and thus $N\cong M[t]$. For a
maximal ideal containing $a$, by choice, $\psi$ is an isomorphism and
thus $N\cong B^n$, in particular extended. Thus $N$ is extended. Also,
notice that $a$ is a non-zero divisor on $N$, since it is free at
maximal ideals containing $a$ and $A$ is an integral domain. Thus, we
have $N_{\mid t=0}=N_0\cong N_1=N_{\mid t=1}$.

Next let us look at $N_0$. Then, we have $M\subset N_0$ and
$N_0/M\cong (A/aA)^n$. Since $N_0$ is a projective module of rank $n$
at primes containing $a$, we see that $M=aN_0$ and since $a$ is a
non-zero divisor on $N_0$, $M\cong N_0$.

Finally, let us look at $N_1$. If we let
$e_1'=e_1-ae_2+a^2e_3-\cdots$, then
$\phi(1)(e_1')=(-1)^{n-1}a^ne_n$. Consider the projection $\pi:A^n\to
A^{n-1}$, to the first $n-1$ factors. Then $\pi\circ\phi(1)$ is onto and
the kernel is generated by $e_1'$. Thus we can identify $N_1$ as the
cokernel of the map 
$$Ae_1'\stackrel{(\phi(1),\psi)}{\longrightarrow} Ae_n\oplus M,$$ 
and 
 $(\phi(1),\psi)(e_1')=((-1)^{n-1}a^ne_n, m_1-am_2+a^2m_3-\cdots)$.
Now, by the previous lemma, we see that $N_1\cong A\oplus
M/(a^n,m)$ and since $N_1\cong N_0\cong M$, we are done.

\medskip
\noindent{\bf Proof  of the Theorem:} With the notation as in the
theorem, we have an inclusion $R\stackrel{(a,m)}{\longrightarrow}
R\oplus M$, which is 
split  and the cokernel is $N$. We wish to show that $M\cong
N$. We will use the following transvections, which do not change the
situation. For any $\phi:M\to R$, we may replace $(a,m)$ by
$(a+\phi(m), m)$. Similarly, for any $m'\in M$, we may replace $(a,m)$
by $(a,m+am')$.

That, $(a,m)$ gives a split inclusion implies  there exists a
homomorphism $\phi:M\to R$ and $b\in R$ such that $ab+\phi(m)=1$. Let
$\mathfrak{m}_1,\ldots, \mathfrak{m}_r$ be the maximal ideals outside
which $M$ is locally free. Let 
$$X=\{\mathfrak{p}\in{\rm Spec\/}\,R\mid a\not\in\mathfrak{p}\,\,
\mbox{and} \,\,\mathfrak{p}\neq \mathfrak{m}_i\}.$$
On this open set $X$, $M$ is locally free and $aM$ generates $M$. So
for a general 
choice of $m'\in M$, $m+am'$ vanishes at only a subset $Z\subset
X$ of codimension $\geq$ rank $M$ by Bertini's theorem.
By assumption 
on the rank of $M$, this 
codimension is at least the dimension of $R$. 
Thus the map $m+am':M^*\to R$ has image height at least $\dim R$
restricted to $X$. 
Also, since $m$ is
unimodular modulo $aM$, $m+am'$ is unimodular at primes
containig $a$ and thus we see that the image of $m+am':M^*\to R$ has
height at least $\dim R$ and comaximal with $a$. We rename $m+am'$ as
$m$ since this is an allowed transvection for our result and thus we
may assume that  
$m(M^*)=I\subset R$ has height at least $\dim R$ and $I$ is comaximal
with $a$. Let
$J=I\cap_i\mathfrak{m}_i$. 
We will arrange $a$ so that $a\not\in \mathfrak{m}_i$
for all $i$. Notice that  for any $\psi\in
M^*$, $a+\psi(m)$ is comaximal with $I$. Let $\phi\in M^*$ be such
that $a$ is comaximal with $\phi(m)$ as before. Though now we have a
new $m$, the same $\phi$ actually works for this $m$ too, though it is
not important. 

We may assume that $a\in\mathfrak{m}_i$ for
$1\leq i\leq p\leq r$ and $a\not\in\mathfrak{m}_i$ for $p<i\leq
r$, possibly after rearranging the $\mathfrak{m}_i$'s. Choose
$x\in\cap_{i>p}\mathfrak{m}_i-\cup_{i=1}^p\mathfrak{m}_i$. Then we may
replace $(a,m)$ by $(a+x\phi(m),m)$. Then we see that $a+x\phi(m)$ is
not in any one of the above maximal ideals. Thus we may assume that
$a$ is comaximal with $J$.

Since characteristic of the field is zero and it is algebraically
closed, by Chinese remainder theorem, we can find a $b\in R$ such that
$a\equiv b^d\pmod J$, where $d$ is the rank of $M$. Since $a=b^d+y$
with $y\in J$ and  since $J\subset I$, there exists a $\psi\in M^*$
such that $\psi(m)=y$. So, we may replace $a$ with $b^d$. Notice that
since $y\in J$ and $a$ is comaximal with $J$, the same holds
for $b$.

Let $K_0$ be the Grothendieck group of projective modules, which is
the same as the Grothendieck group of modules of finite projective
dimension. Thus $[M]\in K_0(R)$ by hypothesis and so we can write
$[M]=[P]-[F]$ for some projective module $P$ over $R$ and $F$, a free
module over $R$. 
Then we have an $x\in J$ so that
$P/(b+x)P$ is free by corollary \ref{coro1}. Let $c=b+x$. Then, 
$c^d=b^d+xz$ for some $z\in R$ and since $x\in J$, there exists a
$\psi\in M^*$ with $\psi(m)=xz$. So we can replace $(b^d,m)$ with
$(c^d,m)$. 
 Since $M$ is projective at all maximal ideals containing $c$, we
see that $M/cM$ is projective and since $c$ is a non-zero divisor in
$M$, we see that $[M/cM]=[P/cP]-[F/cF]$ in $K_0(R/cR)$. By choice of
$c$, $P/cP$ is free over $R/cR$ and thus, the projective module $M/cM$
is stably free over $R/cR$. It has rank
$\geq \dim R>\dim R/cR$. So, 
$M/cM$ is free over $R/cR$ by Bass' cancellation theorem
\cite{Bass3}. Also, the image of $m$ in $M/cM$ is 
unimodular and so by Suslin's theorem \cite{sus2}, we see that the
image of $m$ is part of a free basis of $M/cM$. Now, 
lemma \ref{lem3} finishes the proof.

\section{Reflexive modules over polynomial rings}
Since reflexive modules over polynomial rings in one or two variables
are free (the two variable case was originally proved by
C. S. Seshadri, \cite{CSS}), we will assume that we are over a
polynomial ring in at least three variables. Also, an example in $n$
variables which is not cancellative give an example over polynomial
rings in $k\geq n$ variables, by extending the ring and tensoring the
module.   

\subsection{Four variable case}
Let $R=k[x,y,z,t]$ be a polynomial ring in 4 variables. Let
$v_1=(x,y,zt-1)$ and $v_2=(x,yz,zt-1)$ be two vectors giving rise to
presentations of two modules $M_1,M_2$. Then $M_i$'s are reflexive rank
two modules over $R$. I claim that they are non-isomorphic but
$M_1\oplus R$ is isomorphic to $M_2\oplus R$.

 Notice that $v_1R=v_2R=I$, a complete intersection height
three ideal in $R$. We have exact sequences,
$$0\to M_i^*\to R^3\to I\to 0.$$ We have a commutative diagram,
\[
\begin{array}{ccccccccc}
0&\to &M_1^* &\to& R^3 &\to & I&\to &0\\
&&\uparrow f&&\uparrow \psi &&\uparrow {\rm Id}&&\\
0&\to &M_2^* &\to& R^3 &\to & I&\to &0 
\end{array}
\]
where $\psi$ is the diagonal matrix, $[1,z,1]$. Thus we get an exact
sequence,
$$0\to M_2^*\to M_1^*\to R/zR\to 0.$$ Now by Schanuel's lemma type
argument, we get an exact sequence,
$$0\to R\to R\oplus M_2^*\to M_1^*\to 0.$$
Dualising this and noting that ${\rm Ext\/}^1(M_i^*,R)=R/I$, we see
that,
$$0\to M_1\to M_2\oplus R\to R\to 0$$ is exact and since $R$ is free,
it is split exact. Thus,
$M_1\oplus R\cong M_2\oplus R$.

Now assume that $M_1\cong M_2$. Then we get a commutative diagram,
\[
\begin{array}{ccccccccc}
0&\to & R&\to&R^3&\to&M_1&\to&0\\
&&\downarrow a&&\downarrow \phi&&\downarrow \cong\\
0&\to & R&\to&R^3&\to&M_2&\to&0
\end{array}
\]

With respect to the given bases, we can identify $\phi$ as a $3\times
3$ matrix. 
So, we see that $\det \phi=ca$ where $c\in k$ is a non-zero
constant. Dualising, we see that $a$ is a unit modulo $I$. Now we have
the dual picture as follows:
 
\[
\begin{array}{ccccccccc}
0&\to &M_1^* &\to& R^3 &\to & I&\to &0\\
&&\uparrow \cong&&\uparrow \phi^* &&\uparrow a&&\\
0&\to &M_2^* &\to& R^3 &\to & I&\to &0 
\end{array}
\]

Let us go modulo $I$. Then $I/I^2$ is a free module of rank 3 over
$R/I$ and we have,
\[
\begin{array}{ccc}
(R/I)^3&\to &I/I^2\\
\uparrow\phi^*&&\uparrow a\\
(R/I)^3&\to &I/I^2
\end{array}
\]

We notice that all the maps are now isomorphisms. So, let us compute
$\phi^*$ modulo $I$ with respect to the given  bases. One immediately
sees that $\phi^*$ is the diagonal matrix, $[a,az,a]$. Thus the
determinat of $\phi^*$ is $a^3z$ modulo $I$. But this is equal to
$ca$ since $\det \phi=\det\phi^*$. So, we get $a^3z\equiv ca \bmod I$
or $z\equiv ca^{-2} \bmod 
I$. But since $R/I=k[z,z^{-1}]$, such an equation cannot hold. Thus we
see that $M_1$ and $M_2$ are not isomorphic.

One of the natural questions that can be raised is whether
cancellation does hold if the ranks of the modules are sufficiently
large. But, a modification of the above example shows that it is not
the case. 

For this, take 
$$v_1=(x^n, x^{n-1}y, \ldots, y^n, zt-1)\, \mbox{and} \,  v_2=(zx^n,
x^{n-1}y, \ldots, y^n, zt-1)$$
  and consider the corresponding reflexive
modules $M_1,M_2$. These have rank $n+1$ and as before, it is easy to
show that they are isomorphic, if we add a free module of rank one to
both sides, since the ideal $Rv_1=Rv_2=I$  is of codimension 3 and
$R/I$ is Cohen-Macaulay. If they were isomorphic, exactly as before,
restricting to $R/J$, where $J=(x,y,zt-1)={\rm rad\/}\, I$ we 
get an equation, $ca=a^{n+2}z$ in $R/J=k[z,z^{-1}]$, where $c$ is  a
non-zero constant in $k$ and $a\in R$,  is a unit in $R/J$. 
This leads to a contradiction.

The above examples show that we cannot drop the hypothesis of {\em
locally free outside a finite set of points} in our theorem. 
\subsection{Three variable case}
In three variable polynomial rings, any reflexive module is free
outside a finite set of points. So, if the base fileld is
algebraically closed and the rank is at least three, our theorem will
imply that they are cancellative. Here we construct examples when the
base field is not algebraically closed or if it is of positive
characteristic, of examples which are not cancellative.

\subsubsection{Positive characteristic case}
In this section, we assume that our ring $R=k[x,y,z]$ is a polynomial
ring in three variables and the charactersitic of $k=p>0$. Then, we
will construct examples of reflexive modules of rank $p$ which are
stably isomorphic, but not isomorphic. The method is essentially the
same as above.

Consider, $v_1=(x^{p(p-1)}, x^{p(p-2)}y,\ldots, y^{p-1},z)$ and $v_2$ the
same as $v_1$, except we replace $y^{p-1}$ with $(1+x)y^{p-1}$. Then notice
that $Rv_1=Rv_2=I$ a height three ideal, $R/I$ is Cohen-Macaulay, and
if we denote by $J=(x^p,y,z)$, then $I/JI$ is a free module over $R/J$,
of rank $p+1$. Now the argument is the same as above and if  the
corresponding modules are isomorphic, we get that the element $1+x$
has a $p^{\rm th}$ root in $R/J=k[x]/x^p$, which is
impossible. 
\subsubsection{Characteristic zero case}
Now let us assume that $R=k[x,y,z]$ where $k$ is a field of
characteristic zero and let $n\geq 2$ be an integer. Further assume
that there exists a finite extension $L$ of $k$ such that
$(L^*)^nk^*\neq L^*$ and pick an element $c\in L^*-(L^*)^nk^*$. Using
$c$, as before we will construct  rank $n$ reflexive modules over $R$
which are stably isomorphic, but not isomorphic. Let $p(x)$ be the
irreducible polynomial such that $L=k[x]/p(x)$ and let $q(x)$ be so
chosen so that its image in $L$ via the above map is $c$. Then, we
consider the two vectors,
$$v_1=(p^{n-1},p^{n-2}y,\ldots, y^{n-1}, z), v_2=(p^{n-1},p^{n-2}y,\ldots,
y^{n-1}, qz). $$
One checks as before that the first syzygyies of these two vectors are
rank $n$ reflexive modules, which are stably isomorphic, but not
isomorphic. 

To obtain such examples of fields, we can take for example,
$k=\mathbb{Q}$ and $L=\mathbb{Q}(\theta)$ where
$\theta=\sqrt[n]{2}$. If $(L^*)^nk^*= L^*$, then we get
$\theta=a^nb$ where $a\in L^*$ 
and $b\in k^*$. Taking norms we see that either $2$ or $-2$ is an
$n^{\rm th}$ power in $k$, which is impossible unless $n=1$.

So, we see that the hypothesis of algebraically closed field,
charactersitic zero and locally free outside a finite set of closed
points are all essential for our theorem. Since we have examples
\cite{MK85} of
stably free non-free modules of rank $=\dim R-2$, the condition on the
rank is more or less necessary, except possibly when it is one less
than the dimension of the ring. The condition on finite homological
dimension on the other hand, does not seem essential, but we have
neither a proof nor a counter example.


\begin{thebibliography}{Kum97}

\bibitem[Bas68]{Bass3}
Hyman Bass.
\newblock {\em Algebraic ${K}$-theory}.
\newblock W. A. Benjamin, Inc., New York-Amsterdam, 1968.

\bibitem[Kum85]{MK85}
N.~Mohan Kumar.
\newblock Stably free modules.
\newblock {\em American Journal of Mathematics}, 107:1439--1443, 1985.

\bibitem[Kum97]{MK97}
N.~Mohan Kumar.
\newblock A note on unimodular rows.
\newblock {\em Journal of Algebra}, 191:228--234, 1997.

\bibitem[Qui76]{Q}
D.~Quillen.
\newblock Projective {M}odules over {P}olynomial {R}ings.
\newblock {\em Inventiones {M}athematicae}, 36:166--172, 1976.

\bibitem[Ses58]{CSS}
C.~S. Seshadri.
\newblock Triviality of vector bundles over the affine space $k\sp{2}$.
\newblock {\em Proc. Nat. Acad. Sci. U.S.A.}, 44:456--458, 1958.

\bibitem[Sus76]{sus4}
A.~A. Suslin.
\newblock Projective modules over polynomial rings are free.
\newblock {\em Dokl. Akad. Nauk SSSR}, 229(5):1063--1066, 1976.

\bibitem[Sus77a]{sus1}
A.~A. Suslin.
\newblock A {C}ancellation {T}heorem for {P}rojective {M}odules over
  {A}lgebras.
\newblock {\em Soviet {M}ath. {D}okl.}, 18(5):1281--1284, 1977.

\bibitem[Sus77b]{sus2}
A.A. Suslin.
\newblock Stably free modules.
\newblock {\em Mat. {S}bornik.}, 102(4):537--550, 1977.



\end{thebibliography}
\end{document}